\documentclass[12pt,oneside,english]{amsart}
\usepackage[latin1]{inputenc}
\usepackage{amssymb}

\makeatletter
\newtheorem{theorem}{Theorem}
\newtheorem{lemma}{Lemma}

\newtheorem{definition}{Definition}

\newtheorem{remark}{Remark}
\newtheorem{conjecture}{Conjecture}

\newtheorem{corollary}{Corollary}

\numberwithin{equation}{section}
\usepackage{babel}

\makeatother
\begin{document}
\baselineskip=17pt

\title[ On Stephan's conjectures concerning Pascal triangle]{On Stephan's conjectures concerning Pascal triangle modulo 2 and their polynomial generalization}

\author{Vladimir Shevelev}
\address{Department of Mathematics \\Ben-Gurion University of the
 Negev\\Beer-Sheva 84105, Israel. e-mail:shevelev@bgu.ac.il}


\begin{abstract}
We prove a series of Stephan's conjectures concerning Pascal triangle modulo 2 and give a polynomial generalization.
\end{abstract}

\maketitle

\section{Introduction}
Consider Pascal triangle for binomial coefficient modulo 2. If to read every row of this triangle as a binary number,
then we obtain the following sequence $\{c(n)\}_{n\geq0}$ (cf. A001317 in \cite{11}):
\begin{equation}\label{1.1}
1, 3, 5, 15, 17, 51, 85, 255, 257, 771, 1285, 3855, 4369, 13107, 21845,...
\end{equation}
It is easy to see that
\begin{equation}\label{1.2}
c(2n)\equiv1\pmod4,\enskip n=0,1,...
\end{equation}
Denote
\begin{equation}\label{1.3}
 l(n)=\frac {c(2n)-1} {4}.
 \end{equation}
 In 2004, for sequence $\{l(n)\}_{n\geq0},$ R. Stephan formulated a series of the following conjectures (cf. his comments to A089893 in \cite{11}):
\begin{conjecture}\label{c1}
\begin{equation}\label{1.4}
l(2^k) = 2^{2^{k+1}-2}.
\end{equation}
\end{conjecture}
\begin{conjecture}\label{c2}
\begin{equation}\label{1.5}
\lim_{n\rightarrow\infty} l(2n+1)/l(2n)) = 5.
\end{equation}
\end{conjecture}
\begin{conjecture}\label{c3}
\begin{equation}\label{1.6}
\lim_{n\rightarrow\infty} l(4n+2)/l(4n+1)) = 17/5.
\end{equation}
\end{conjecture}
\begin{conjecture}\label{c4}
\begin{equation}\label{1.7}
\lim_{n\rightarrow\infty} l(8n+4)/l(8n+3)) = 257/85.
\end{equation}
\end{conjecture}
etc.\newline
\indent We add that Moscow PhD student S. Shakirov conjectured (private communication) that a generating function for sequence $\{c(n)\}$ is
\begin{equation}\label{1.8}
\prod_{k=0}^{\infty}(1+x^{2^k}+(2x)^{2^k})=\sum_{n=0}^{\infty}c(n)x^n.
\end{equation}

\indent In this paper we prove these conjectures and give a polynomial generalizations.
   \section{On sequence A001317 }
    Consider an infinite in both sides $(0,1)$-sequence with a finite set of 1's which we call $C$-sequence. Removing in it all 0's before the first 1 and after the last 1, we obtain some odd number which we call the kernel of  $C$-sequence. Every $C$-sequence generates a new $C$-sequence, if to write sums of every pair of its adjacent terms modulo 2. If to consider infinite iterations of such process beginning with $C$-sequence with kern 1, then we obtain $C$-sequences, the kernels $\{c(i)\}_{i\geq0}$ of which form Pascal's triangle for binomial coefficients modulo 2.  Note that, $c(0)=1$ and $c(i)$ contains $i+1$ binary digits.\newline
    \indent Consider now sequence $\{d(n)\}$ defined by the formula $d(0)=1;$ for $n\geq1,$ if binary expansion of $n$ is
     \begin{equation}\label{2.1}
    n=\sum_{i=1}^{m} 2^{k_i},
    \end{equation}
    then
    \begin{equation}\label{2.2}
 d(n)=\prod_{i=1}^{m}F(k_i),
\end{equation}
where
\begin{equation}\label{2.3}
F(n)=2^{2^n}+1, \enskip n\geq0,
\end{equation}
is Fermat number. Such decomposition of $d(n)$ we call its \slshape Fermat factorization. \upshape \newline
\indent From (\ref{2.1})-(\ref{2.2}) immediately follows a generating function for $\{d(i)\}:$
\begin{equation}\label{2.4}
\prod_{k=0}^{\infty}(1+F(k)x^{2^k})=\sum_{n=0}^{\infty}d(n)x^n,\enskip 0<x<\frac{1} {2}.
\end{equation}
Note that sequence $\{d(i)\}$ possesses the following properties:\newline
1) $d(n)$ is a binary number with $n+1\enskip (0,1)$-digits;\newline
2) numbers $\{d(i)\}$ are 1 and  all Fermat numbers or products of distinct Fermat numbers;\newpage
3) number of Fermat factors in the product equals to $d(n)$ is the number of 1's in the binary expansion of $n.$\newline
4) $F(i)$ divides $d(n),\enskip n>1,$ if and only if it is a factor in product (\ref{2.2}). \newline
\indent Proofs of these properties is very easy: 1) follows from a simple induction; 2) and 3) follow from the definition; 4) follows from the well known fact (cf., e.g., \cite{12}) that every two Fermat numbers are relatively prime, in view of recursion
\begin{equation}\label{2.5}
F(n)=2+\prod_{i=0}^{n-1}F(i).
\end{equation}
\begin{theorem}\label{t1} For $n=0,1,...,$ we have
\begin{equation}\label{2.6}
c(n)=d(n).
\end{equation}
\end{theorem}
\bfseries Proof. \mdseries We use induction, the base of which is $c(0)=d(0)=1, \enskip c(1)=d(1)=3,\enskip c(2)=d_2=5.$ Suppose that $c(i)=d(i),$ for $i\leq k.$ Let $m$ be the most number for which $F(m)$ divides $c(k)=d(k).$ In non-trivial case, when $c(k)\neq F(m),$ using property 4), for some $r<k,$ we have $c(k)=d(r)F(m)=c(r)F(m).$ Furthermore, since, by the condition, $F(m)$ is the most Fermat divisor of $c(k)$ and, in view of (\ref{2.5}), we have
\begin{equation}\label{2.7}
 c(r)=\frac {c(k)}{F(m)}\leq \prod_{i=0}^{m-1}F(i)=F(m)-2.
 \end{equation}
Besides, since $c(r)<c(k),$ then, by the inductive supposition,
$$c(r+1)=d(r+1).$$
Adding the case when $c(k)=F(m),$ let us prove a recursion: $c(0)=1, c(1)=3, c(2)=5;$ for $k\geq2,$
\begin{equation}\label{2.8}
c(k+1)=\begin{cases}3F(m),\enskip if\enskip c(k)=F(m),\\ F(m+1), \enskip if \enskip 1<c(r)=F(m)-2,\\F(m)c(r+1), \enskip if \enskip 1<c(r)<F(m)-2.\end{cases}
\end{equation}
\indent Let $c(k)=F(m),\enskip m\geq1. \enskip C$-sequence with kernel $c(k)$ is
$$...01\underbrace{0...0}_{2^m-1}10...$$
Thus the following $C$-sequence with kernel $c(k+1)$ is
$$...011\underbrace{0...0}_{2^m-2}110...$$
Comparing kernels $c(k)$ and $c(k+1),$ we conclude that $c(k+1)=3c(k)=3F(m).$\newpage
\indent Furthermore, if $c(r)=F(m)-2,$
then, by (\ref{2.7}), we have
$$c(k)=F(m)c(r)=F(m)(F(m)-2)=F(m+1)-2=\underbrace{11...1}_{2^{m+1}}.$$
Thus the $C$-sequence with kernel $c(k)$ is
$$...0\underbrace{11...1}_{2^{m+1}}0...$$
Therefore, by the definition, the $C$-sequence with kernel $c(k+1)$ is
$$...01\underbrace{0...0}_{2^{m+1}-1}10...$$
and we see that $c(k+1)=F(m+1).$ \newline
\indent Let now $c(r)<F(m)-2.$ Since, by the supposition of induction,  $c(r)=d(r).$ Therefore, $c(r)$ is a product of Fermat numbers and
$$c(r)\leq \frac{\prod_{i=0}^{m-1}F(i)} {F(0)}=\frac{F(m)-2}{F(0)}.$$
Hence, $c(r)$ is not more than $(2^m-1)$-digits odd binary number.
Since
$$c(k)=F(m)c(r)=2^{2^m}c(r)+c(r),$$
then $c(k)$ has the binary expansion of the form
\begin{equation}\label{2.9}
c(k)=\overline{c(r)\underbrace{0...0}_lc(r)},
\end{equation}
where $l\geq1.$ \newline
Passing on to the following kernel, we have:
$$c(k+1)=\overline{c(r+1)\underbrace{0...0}_{l-1}c(r+1)},$$
where $l-1\geq0.$ Thus
$$c(k+1)=c(r+1)2^{2^m}+c(r+1)=c(r+1)F(m). $$
This completes formula (\ref{2.8}). From this formula we conclude that $c(k+1)$ is a term of sequence $\{d(i)\}.$
Moreover, since $c(k+1)$ contains $k+2$ binary digits, then, in view of property 1) of numbers $\{d(i)\},$ both of $c(k+1)$ and $d(k+1)$ contain $(k+2)$ binary digits. Therefore, $c(k+1)=d(k+1). \blacksquare$
\begin{remark}\label{r1}
In proof of Theorem \ref{t1} we essentially followed to our arguments from preprint \cite{9}, $1991.$
\end{remark}
\begin{remark}\label{r2}
Hewgill \cite{4}, for the first time, found a relationship between Pascal's triangle modulo $2$ and Fermat numbers. In fact, using a simple induction, he proved the following explicit formula for the binary representation of $c_n$:
\newpage
$$c_n=(\prod_{i=0}^{\lfloor\log_2n\rfloor}F_i^{(\lfloor\frac{n}{2^i}\rfloor \pmod {2})} )_2.$$
\end{remark}
\begin{remark}\label{r3}
Karttunen \cite{6} gave a representation of $c_n$ in the Fibonacci number system.
\end{remark}
\begin{corollary}\label{cor1}
Conjectural generating formula $(\ref{1.8})$ is true.
\end{corollary}
\bfseries Proof. \mdseries According to (\ref{2.4}) and Theorem \ref{t1}, we have
\begin{equation}\label{2.10}
\prod_{k=0}^{\infty}(1+F(k)x^{2^k})=\sum_{n=0}^{\infty}c(n)x^n,\enskip 0<x<\frac{1} {2}.
 \end{equation}
 It is left to note that
 $$1+F(k)x^{2^k}=1+x^{2^k}+(2x)^{2^k}. \blacksquare $$
\indent Denote $s(n)$ the number of 1's
in the binary expansion of $n.$
\begin{corollary}\label{cor2}
$a)$ Number of factors in Fermat factorization of $c(n)$ is $s(n).$\newline
$b)$ Moreover, the following formula holds
\begin{equation}\label{2.11}
 s(c(n))=2^{s(n)}.
\end{equation}
\end{corollary}
\bfseries Proof. \mdseries a) follows from Theorem \ref{t1} and property 3) of numbers $\{d(n)\}.$\newline
\indent b) Let, firstly, $c(k)$ be not a Fermat number and, as in proof of Theorem \ref{t1}, $m$ be the most number for which $F(m)$ divides $c(k),$ such that $c(k)=F(m)c(r).$ Since the difference between numbers of factors in Fermat factorization of $c(k)$ and $c(r)$ is 1, then, according to a), we have
$$s(k)=s(r)+1.$$
Now we use induction. If the statement is true for $i\leq k-1,$ then, in particular, $ s(c(r))=2^{s(r)}.$ Therefore, by (\ref{2.9}), we have
$$s(c(k))=2s(c(r))=2\cdot2^{s(r)}=2^{s(r)+1}=2^{s(k)}.$$
It is left to consider case $c(k)=F(l).$ Here, by a), $s(k)=1$ and (\ref{2.9}) satisfies trivially. $\blacksquare$\newline
\indent Note that point b) of Corollary \ref{cor2} means that the number of odd binomial coefficient in $n$-th row of Pascal triangle is $2^{s(n)}.$ It is known result of J.Glaisher \cite{2}. His proof was based on well known Lucas (1878) comparison modulo 2:
if the binary representations of numbers $m\geq t$ are $m=\overline{m_1...m_k}, \enskip t=\overline{t_1...t_k}$ (with, probably, some first $t_i=0$), then

$$\binom {n} {t}\equiv\prod_{i=0}^{m}\binom {n_i} {t_i}\pmod 2.$$
\newpage
In \cite{3} A.Granville gives a new interesting proof of Glaisher's result. Our proof is the third one. Generalizations in other directs see in \cite{1}, \cite {3}, \cite{5}, \cite{8}, \cite{10}.
\begin{corollary}\label{cor3}
If $F(m)$ is the most Fermat divisor of numbers $c(k-1)$ and $c(l-1)$ from interval $(1,\enskip F(m)-2),$ then
\begin{equation}\label{2.12}
c(k-1)c(l)=c(l-1)c(k).
\end{equation}
\end{corollary}
\bfseries Proof. \mdseries Using (\ref{2.8}), we have
$$c(k)=c(k-1)F(m),\enskip c(l)=c(l-1)F(m)$$
and (\ref{2.12}) follows. $\blacksquare$
\begin{corollary}\label{cor4}
If $k=2^ml+2^{m-1},\enskip m\geq1,$ then
\begin{equation}\label{2.13}
c(k)=c(2^ml)F(m-1).
\end{equation}
\end{corollary}
\bfseries Proof. \mdseries From (\ref{2.1})-(\ref{2.2}), we immediately have $d(k)=d(2^ml)F(m-1),$ and (\ref{2.13}) follows from Theorem \ref{t1}. $\blacksquare$

\section{Proof of Conjecture \ref{c1}}
Now proof of Conjecture \ref{c1} is especially simple. Indeed, in view of (\ref{1.3}) and (\ref{2.3}), formula (\ref{1.4}) of Conjecture \ref{c1} can be rewritten as
 \begin{equation}\label{3.1}
 c(2^n)=F(n),
\end{equation}
where $n=k+1\geq1.$\newline
\indent According to Corollary \ref{cor1}$a),$ number $c(2^n)$ has only one Fermat factor, i.e., for some $t,$ we have
$c(2^n)=F_t.$ Besides, by the definition, $c(2^n)$ has $2^n+1$ binary digits. It is left to notice that, the unique Fermat number having $2^n+1$ binary digits is $F(n),$ i.e., $t=n$ and $c(2^n)=F(n). \blacksquare$\newline
\indent In addition, prove that
\begin{equation}\label{3.2}
c(2^n-1)=F(n)-2.
\end{equation}
Indeed, by the definition of sequence $\{d(n)\}$ and (\ref{2.3}), we conclude that $F(n)-2,$ as a product of \slshape distinct \upshape Fermat numbers, is a term of sequence $\{d(i)\}$ and thus, by Theorem \ref{t1}, is a term of sequence $\{c(i)\}.$  Now it is left to notice that numbers $c(2^n-1)$ and $F(n)-2$ have the same number $(2^n)$ of binary digits. $\blacksquare$

\section{Proof of Conjectures \ref{c2}, \ref{c3}, \ref{c4}, $etc.$}

\begin{lemma}\label{l1} For every $n\geq0,\enskip t\geq1$ we have identity
\begin{equation}\label{4.1}
(F(t-1)-2)c(2^{t}n)=c(2^{t}n+2^{t-1}-1).
\end{equation}
\end{lemma}
\newpage
\bfseries Proof. \mdseries As in proof of (\ref{3.2}), we conclude that $(F(t-1)-2)c(2^{t}n)$ is a term of sequence $\{c(i)\}.$ Note that number $c(2^{t}n+2^{t-1}-1)$ has $2^{t}n+2^{t-1}$ binary digits. Besides, number $F(t-1)-2=\underbrace{1...1}_{2^{t-1}}$ and $c(2^{t}n)$ has $2^{t}n+1$ binary digits. Therefore, number $(F(t-1)-2)c(2^{t}n)$ contains not less binary digits than number $\underbrace{1...1}_{2^{t-1}}\underbrace{0...0}_{2^{t}n},$ i.e. $(F(t-1)-2)c(2^{t}n)$ has not less than $2^{t-1}+2^{t}n$ binary digits. On the other hand, $(F(t-1)-2)c(2^{t}n)$ contains not more binary digits than number  $$\underbrace{1...1}_{2^{t-1}}\underbrace{1...1}_{2^{t}n}=(2^{2^{t-1}}-1)(2^{2^tn}-1)\leq 2^{2^{t-1}+2^tn}-1, $$
 i.e., $(F(t-1)-2)c(2^{t}n)$ has not more than $2^{t-1}+2^{t}n$ binary digits. Thus number $(F(t-1)-2)c(2^{t}n)$ has
 exactly $2^{t-1}+2^{t}n$ binary digits. Consequently, two terms $(F(t-1)-2)c(2^{t}n)$ and $c(2^{t}n+2^{t-1}-1)$ of sequence $\{c(i)\}$ has the same number of digits. Therefore, equality (\ref{4.1}) holds. $\blacksquare$
 \begin{lemma}\label{l2} For every $n\geq0,\enskip t\geq1,$ we have identities
\begin{equation}\label{4.2}
(F(t-1)-2)c(2^{t}n+2^{t-1})=F(t-1)c(2^{t}n+2^{t-1}-1),
\end{equation}
\begin{equation}\label{4.3}
(F(t-1)-2)c(2^{t}n+2^{t-1})=3F(t-1)c(2^{t}n+2^{t-1}-2).
\end{equation}
\end{lemma}
\bfseries Proof. \mdseries Multiplying (\ref{4.1}) by $F(t-1)$ and using formula (\ref{2.13}) of Corollary \ref{cor4} (for $l=n$ and $m=t$), we obtain (\ref{4.2}). Furthermore, if to take in Corollary \ref{cor4} $m=1,\enskip l=2^{t-1}n+2^{t-2}-1,$ then, in view of $F(0)=3,$ we have $c(2^{t}n+2^{t-1}-1)=3c(2^{t}n+2^{t-1}-2),$ and (\ref{4.3}) follows. $\blacksquare$ \newline\newline

\indent Now we are able to get a proof of Conjectures \ref{c2}, \ref{c3}, \ref{c4}, $etc.$ According to (\ref{1.3}), we have
\begin{equation}\label{4.4}
c(2n)=4l(n)+1.
\end{equation}
Let in (\ref{4.3}) $t\geq2.$ Then, by (\ref{4.4}), we have
$$(F(t-1)-2)(4l(2^{t-1}n+2^{t-2})+1)=3F(t-1)(4l(2^{t-1}n+2^{t-2}-1)+1),$$
or
$$ \frac {4l(2^{t-1}n+2^{t-2})+1}{4l(2^{t-1}n+2^{t-2}-1)+1}=\frac{3F(t-1)}{F(t-1)-2}. $$
Hence, we finally find
\begin{equation}\label{4.5}
\lim_{n\rightarrow\infty} \frac{l(2^{t-1}n+2^{t-2})}{l(2^{t-1}n+2^{t-2}-1)}=\frac{3F(t-1)}{F(t-1)-2}.
\end{equation}
 $\blacksquare$\newpage
\indent So, if $t=2,3,4,5,...,$ then the right hand side is
$$\frac {3\cdot5 } {5-2 }=5, \enskip \frac {3\cdot{17}} {17-2 }=\frac{17} {5},\enskip \frac {3\cdot{257}} {257-2 }=\frac{257} {85},\enskip \frac {3\cdot{65537}} {65537-2 }=\frac{65537} {21845}, ...$$ respectively.

\section{Second proof of key identity (\ref{4.3}) based on notion of orthogonality of nonnegative integers}
We can essentially simplify our proof of Stephan's conjectures by a simplification of key identity (\ref{4.3}).
Put to every nonnegative integer $n$ to one-to-one correspondence $(0,1)$-vector $\overline{n}$ by the rule: if the binary expansion of $n$ is $n=\overline{n_1...n_m},$ then
\begin{equation}\label{5.1}
\overline{n}=\overline{...0...0n_1...n_m}
\end{equation}
with infinitive 0's before $n_1.$
For two integers $u\leq v$ with vectors $\overline{u}=\overline{...0...0u_1...u_l}$ and $\overline{v}=\overline{...0...0v_1...v_m},\enskip l\leq m$ introduce "circ-product" by formula ( which is, for the corresponding vectors, similar to dot-product)
\begin{equation}\label{5.2}
u\circ v=\overline{u}\overline{v}=u_lv_m+u_{l-1}v_{m-1}+...+u_1v_{m-l+1}.
\end{equation}
\begin{definition}\label{d1}
We call two non-negative integers \enskip$u,\enskip v$ \upshape mutually orthogonal $(u\bot v),$ \slshape if
$u\circ v=0.$
\end{definition}
Note that if $(u\bot v),$ then the sets of positions of 1's in their binary representations do not intersect.

An important source for obtaining various identities for numbers $\{c(n)\}$ is the following exponential-like "addition theorem".
\begin{lemma}\label{l3}
If $n_1\bot n_2,$ then
\begin{equation}\label{5.3}
c(n_1+n_2)=c(n_1)c(n_2).
\end{equation}
\end{lemma}
\bfseries Proof. \mdseries Let $n_1\geq n_2$ and the binary expansions of $n_1$ and $n_2$ be $n_1=\sum_{i=1}^{m} 2^{k_i}$ and $n_2=\sum_{j=1}^{m} 2^{l_j}$  (with, probably, some first $l_i=0$). Since $n_1\bot n_2,$ then $k_i\neq l_j,\enskip i,j=1,...,m.$ Thus the binary expansion of $n_1+n_2$ is $\sum_{i=1}^{m}2^{k_i}+\sum_{j=1}^{m}2^{l_j}.$
Therefore, according to (\ref{2.1})-(\ref{2.2}), we have
$$c(n_1+n_2)=(\prod_{i=1}^{m}F(k_i))(\prod_{j=1}^{m}F(l_j))=c(n_1)c(n_2). \blacksquare$$
\bfseries Second proof of (\ref{4.3}). \mdseries \newline
a)Using the notion of numbers orthogonality, we immediately obtain formula (\ref{4.2}) by the following way.\newline
 By (\ref{3.2}), we have
 \newpage
\begin{equation}\label{5.4}
F(t-1)-2=c(2^{t-1}-1).
\end{equation}
Since, evidently, $(2^{t-1}-1)\bot(2^tn+2^{t-1}),$ then, using (5.3)-(5.4), we find
$$(F(t-1)-2)c(2^{t}n+2^{t-1})=c(2^{t}n+2^{t-1}+2^{t-1}-1)=c(2^tn+2^t-1).$$
On the other hand, since $2^{t-1}\bot(2^tn+2^{t-1}-1),$ then
$$F(t-1)c(2^{t}n+2^{t-1}-1)=c(2^{t-1})c(2^{t}n+2^{t-1}-1)=c(2^tn+2^t-1).$$
Thus we conclude that (\ref{4.2}) holds. \newline
 b) Note now that, $1\bot 2^tn+2^t-2.$ Thus $3c(2^{t}n+2^{t-1}-2)=c(2^{t}n+2^{t-1}-1)$ and (\ref{4.3}) follows as well.$\blacksquare$\newline
\indent Further we consider a polynomial generalization.

\section{Polynomials $p_n(z),\enskip q_n(z)$ and their properties }
Consider sequence of polynomials (cf.\cite{3})
\begin{equation}\label{6.1}
p_n(z)=\frac {1} {2}\sum_{i=0}^n(1-(-1)^{\binom {n} {i}})z^i, \enskip n=0,1,..., \enskip z\in\mathbb{C},
\end{equation}
such that
\begin{equation}\label{6.2}
p_n(0)=1,\enskip p_n(1)=2^{s(n)},   \enskip p_n(2)=c(n).
\end{equation}
The second equality we have in view of (\ref{2.9}).\newline
\indent By the same way, one can prove a generalization of Theorem \ref{t1}.
\begin{theorem}\label{t2} For $n\geq1,$ we have the following decomposition of $p_n(z):$
\begin{equation}\label{6.3}
p_n(z)=\prod_{i=0}^m(z^{2^{k_i}}+1),
\end{equation}
if the binary expansion of $n$ is
\begin{equation}\label{6.4}
n=\sum_{i=0}^m2^{k_i}.
\end{equation}
Thus a generating function for polynomials $\{p_n(z)\}$ is
\begin{equation}\label{6.5}
\prod_{k=0}^{\infty}(1+(z^{2^{k}}+1)x^{2^k})=\sum_{n=0}^{\infty}p_n(z)x^n,\enskip 0<x<\frac{1} {|z|}.
\end{equation}
\end{theorem}
In particular, we have
\begin{equation}\label{6.6}
p_{2^n}(z)=z^{2^{n}}+1.
\end{equation}
Note that, if $n$ has binary expansion (\ref{6.4}), then $2n=\sum_{i=0}^m2^{k_i+1}.$ Since $z^{2^{k_i+1}}=(z^2)^{2^{k_i}},$ then we have
\begin{equation}\label{6.7}
p_{2n}(z)=\prod_{i=0}^m((z^2)^{2^{k_i}}+1)=p_n(z^2).
\end{equation}
\newpage
Analogously, since $2n+1=1+\sum_{i=0}^m2^{k_i+1},$ then
\begin{equation}\label{6.8}
p_{2n+1}(z)=(z+1)\prod_{i=0}^m((z^2)^{2^{k_i}}+1)=(z+1)p_n(z^2).
\end{equation}
Formulas (\ref{6.7})-(\ref{6.8}) give a simple recursion for polynomials $\{p_n(z)\},$ which recently were obtained by
S. Northshield (cf. \cite{7}, Lemma 3.1) in a quite another way.\newline
\indent Note that every two different polynomials in sequence $\{p_{2^{i}}(z)=z^{2^{i}}+1\}_{i\geq0}$ are respectively prime. It follows from the identity
\begin{equation}\label{6.9}
p_{2^{n}}(z)=2+(z-1)\prod_{i=0}^{n-1}p_{2^{i}}(z).
\end{equation}
Put
\begin{equation}\label{6.10}
F_n(z)=p_{2^n}(z)=z^{2^{n}}+1.
\end{equation}
 The following identity holds (cf. \cite{9})
\begin{equation}\label{6.11}
\sum_{n=0}^{\infty}\frac {1} {p_n(z)^s}=\prod_{k=0}^\infty (1+F_k(z)^{-s}), \enskip |z|>1, \enskip \Re{s}>0.
\end{equation}
In particular, for $z=2,\enskip s=1,$ we have
\begin{equation}\label{6.12}
\sum_{n=0}^{\infty}\frac {1} {c(n)}=\prod_{k=0}^\infty (1+F_k^{-1})=1.700735495...\enskip .
\end{equation}
According to Theorem \ref{t2} and in view that $s(n)\equiv m_n \pmod2,$ where ${m_n}={0,1,1,0,1,0,0,1,1...}$ is Thou-Morse sequence, together with (\ref{6.11}), we have also
\begin{equation}\label{6.13}
\sum_{n=0}^{\infty}\frac {(-1)^{m_n}} {p_n(z)^s}=\prod_{k=0}^\infty (1-F_k(z)^{-s}), \enskip |z|>1, \enskip \Re{s}>0.
\end{equation}
Let us show that, in particular, for $s=1,$ we have
\begin{equation}\label{6.14}
\sum_{n=0}^{\infty}\frac {(-1)^{m_n}} {p_n(z)}=1-\frac{1} {z}, \enskip |z|>1.
\end{equation}
Indeed, since
$$1-\frac{1} {F_n(z)}=(1+\frac{1} {z^{2^n}})^{-1}, $$
then
$$\prod_{k=0}^\infty (1-F_k(z)^{-1})=\prod_{k=0}^\infty(1+\frac{1}{z^{2^n}})^{-1} $$
and it is left to note that
\begin{equation}\label{6.15}
\prod_{n=0}^\infty (1+\frac{1} {z^{2^n}})=1-\frac{1} {z}.
\end{equation}
\newpage
In particular, together with (\ref{6.12}), for $z=2,$ we find
\begin{equation}\label{6.16}
\sum_{n=0}^{\infty}\frac {(-1)^{m_n}} {c(n)}=\frac{1} {2}.
\end{equation}
In addition, note that, if to consider all different finite products of not necessarily distinct polynomials from sequence $\{p_n(z)\},$ then we obtain a sequence of polynomials $q_n(z):$
$$q_0(z)=1,\enskip q_1(z)=z+1,\enskip q_2(z)=z^2+1,\enskip q_3(z)=(z+1)^2, $$
\begin{equation}\label{6.17}
  q_4(z)=(z+1)(z^2+1),\enskip q_5(z)=z^4+1,\enskip q_6(z)=(z^2+1)^2.
 \end{equation}
 For these polynomials, together with (\ref{6.11}), we have the following analog of Euler identity for primes:
 \begin{equation}\label{6.18}
 \prod_{F\in F(z)}(1-F^{-s})^{-1}=\sum_{n=0}^{\infty}\frac {1}{q_n(z))^s},\enskip |z|>1, \enskip \Re{s}>0,
 \end{equation}
 where
  $$F(z)=\{F_n(z)\}_{n\geq0}.$$
 In particular, for $s=1,$ using (\ref{6.15}), we have
$$\sum_{n=0}^{\infty}\frac {1}{q_n(z)}=\prod_{F\in F(z)}(1-F^{-1})^{-1}=$$
 \begin{equation}\label{6.19}
\prod_{n=0}^\infty (1+\frac{1} {z^{2^n}})^{-1}=\frac{z} {z-1},\enskip |z|>1.
  \end{equation}
 Furthermore, introducing an analog of M\"{o}bius function
  \begin{equation}\label{6.20}
\nu (n)=\begin{cases}
(-1)^{m_n},\;\; if\;\;n\enskip is\enskip squarefree,\\0,\;\;otherwise,\end{cases}
\end{equation}
we get
 \begin{equation}\label{6.21}
\sum_{n=0}^{\infty}\frac {\nu(n)} {q_n(z)^s}=\prod_{F\in F(z)}(1-F^{-s}),\enskip |z|>1, \enskip \Re{s}>0.
\end{equation}
In particular, for $s=1,$ we have

 \begin{equation}\label{6.22}
\sum_{n=0}^{\infty}\frac {\nu(n)}{q_n(z)}=1-\frac{1} {z},\enskip |z|>1.
  \end{equation}
  \newpage
\section{Polynomial generalization of Stephan's relations}
Now we consider a polynomial generalization of formulas of the previous sections which leads us to the corresponding generalization of Stephan's relations. Since proof of the generalized formulas is quite analogous, then we restrict ourself only by writing of the chain of them. For $|z|>1,$ we have
\begin{equation}\label{7.1}
p_{2^n-1}=\frac{F_n(z)-2} {z-1}.
\end{equation}
This formula generalizes (\ref{3.2}). Furthermore, the following generalization of (\ref{2.13}) holds:
\begin{equation}\label{7.2}
p_{2^ml+2^{m-1}}(z)=p_{2^ml}(z)F_{m-1}(z).
\end{equation}
In particular, taking in (\ref{7.2}) $m=1,\enskip l=2^{t-1}n+2^{t-2}-1,$ in view of $F_0(z)=z+1,$ we find
\begin{equation}\label{7.3}
p_{2^tn+2^{t-1}-1}(z)=(z+1)p_{2^tn+2^{t-1}-2}(z).
\end{equation}
After that the corresponding generalization of formulas (\ref{4.1})-(\ref{4.3}) is obtained. We have
\begin{equation}\label{7.4}
(F_{t-1}(z)-2)p_{2^tn}(z)=p_{2^tn+2^{t-1}-1}(z),
\end{equation}
\begin{equation}\label{7.5}
(F_{t-1}(z)-2)p_{2^{t}n+2^{t-1}}(z)=(z-1)F_{t-1}(z)p_{2^tn+2^{t-1}-1}(z),
\end{equation}
\begin{equation}\label{7.6}
(F_{t-1}(z)-2)p_{2^{t}n+2^{t-1}}(z)=(z^2-1)F_{t-1}(z)p_{2^tn+2^{t-1}-2}(z).
\end{equation}
Note that
\begin{equation}\label{7.7}
p_{2n}(z)\equiv1\pmod{z^2}.
\end{equation}
Put
\begin{equation}\label{7.8}
l_n(z)=\frac{p_{2n}(z)-1}{z^2}.
\end{equation}
Let in (\ref{7.6}) $t\geq2.$
Then we have
\begin{equation}\label{7.9}
(F_{t-1}(z)-2)(z^2l_{2^{t-1}n+2^{t-2}}(z)+1)=(z^2-1)F_{t-1}(z)(z^2l_{2^{t-1}n+2^{t-2}-1}(z)+1),
\end{equation}
or
\begin{equation}\label{7.10}
\frac {z^2l_{2^{t-1}n+2^{t-2}}(z)+1} {z^2l_{2^{t-1}n+2^{t-2}-1}(z)+1}=\frac{(z^2-1)F_{t-1}(z)}{F_{t-1}(z)-2}
\end{equation}
and, consequently,
\begin{equation}\label{7.11}
\lim_{n\rightarrow\infty}\frac {l_{2^{t-1}n+2^{t-2}}(z)} {l_{2^{t-1}n+2^{t-2}-1}(z)}=\frac{(z^2-1)F_{t-1}(z)}{F_{t-1}(z)-2}.
\end{equation}

In particular, for $t=2,$
\newpage
$$\lim_{n\rightarrow\infty}\frac {l_{2n+1}(z)} {l_{2n}(z)}=z^2+1;$$
for $t=3,$
$$\lim_{n\rightarrow\infty}\frac {l_{4n+2}(z)} {l_{4n+1}(z)}=\frac{z^4+1}{z^2+1};$$
for $t=4,$
$$\lim_{n\rightarrow\infty}\frac {l_{8n+4}(z)} {l_{8n+3}(z)}=\frac{z^8+1}{(z^4+1)(z^2+1)},\enskip etc.$$

In case of $z=2,$ we again obtain formulas (\ref{1.5})-(\ref{1.7}).

\;\;\;\;\;\;\;\;
\end{document}